\documentclass[twoside]{amsart}
\usepackage{latexsym}
\usepackage{amssymb,amsmath,amsopn}
\usepackage[dvips]{graphicx}   
\usepackage{color,epsfig}      

\newtheorem*{thm}{Theorem}
\newtheorem*{lem}{Lemma}
\theoremstyle{definition}
\newtheorem*{defn}{Definition} 
\newtheorem*{rem}{Remark}
\newtheorem{ques}{Question}
\newtheorem{prob}{Problem}
\newtheorem*{example}{Example}

\newcommand{\Eu}{\mathbb E}
\newcommand{\M}{\mathbb{M}}
\renewcommand{\S}{\mathbb{S}}
\newcommand{\HH}{\mathbb{H}}

\DeclareMathOperator{\inter}{int}
\DeclareMathOperator{\bd}{bd}
\DeclareMathOperator{\dist}{dist}

\DeclareMathOperator{\perim}{perim}
\DeclareMathOperator{\area}{area}
\DeclareMathOperator{\diam}{diam}

\begin{document}
\title[On the perimeters of simple polygons]{On the perimeters of simple polygons contained in a disk}

\author[Z. L\'angi]{Zsolt L\'angi}
                                                                               
\address{Zsolt L\'angi, Dept. of Geometry, Budapest University of Technology,
Budapest, Egry J\'ozsef u. 1., Hungary, 1111}
\email{zlangi@math.bme.hu}
                                                                        

\begin{abstract}
A simple $n$-gon is a polygon with $n$ edges with each vertex belonging to exactly two edges and
every other point belonging to at most one edge.
Brass \cite{BMP05} asked the following question: For $n \geq 5$ odd,
what is the maximum perimeter of a simple $n$-gon contained in a Euclidean unit disk?
In 2009, Audet, Hansen and Messine \cite{AHM09} answered this question, and showed that the optimal configuration
is an isosceles triangle with a multiple edge, inscribed in the disk.
In this note we give a shorter and simpler proof of their result, which
we generalize also for hyperbolic disks, and for spherical disks of sufficiently small radii.
\end{abstract}
\maketitle

\section{Introduction}

This paper deals with questions about the Euclidean plane, denoted by $\Eu^2$,
the hyperbolic plane, denoted by $\HH^2$, and the sphere, denoted by $\S^2$.

A question regarding an isoperimetric problem about simple polygons
was asked by Peter Brass (see Problem 3 on p. 437 in \cite{BMP05}).

\begin{prob}[Brass, 2005]
For $n \geq 5$ odd, what is the maximum perimeter of a simple $n$-gon contained in a Euclidean unit disk?
\end{prob}

The authors of \cite{BMP05} remarked that for $n$ even, the supremum of the perimeters is
the trivial upper bound $2n$, as it can be approached by simple $n$-gons in which the sides go
back and forth near a diameter of the disk.
This argument does not work for an odd value of $n$.
Nevertheless, it is easy to show that the maximum perimeter of a triangle is attained for the regular triangle.
In 2009, Audet, Hansen and Messine \cite{AHM09} showed that for $n$ odd, the optimal configuration is an isosceles
triangle with a multiple edge, inscribed in a disk.

In this paper, we present a shorter and simpler proof of their result, and generalize it for a hyperbolic disk of arbitrary radius, and also for a spherical disk of sufficiently small radius.

We use the following notations.
Let $\M \in \{ \Eu^d, \HH^d, \S^d \}$, and let $x, y\in \M$.
The distance of $x$ and $y$ is denoted by $\dist_M(x,y)$.
If $\M \neq \S^d$, or $x$ and $y$ are not antipodal points of $\S^d$, then
$[x,y]$ (respectively, $(x,y)$) denotes the closed (respectively, open) segment with endpoints $x$ and $y$,
and if $x \neq y$, $L(x,y)$ denotes the straight line passing through $x$ and $y$.

For any set $A \subset \M$, we use the standard notations
$\inter A$, $\bd A$, $\diam A$, $\perim A$ and $\area A$ for the \emph{interior}, the \emph{boundary},
the \emph{diameter}, the \emph{perimeter} or the \emph{area} of $A$.
Points are denoted by small Latin letters, sets of points by capital Latin letters, and
real numbers by Greek letters.

Our main results are the following.

\begin{thm}\label{thm:boundary}
Let $\M$ be $\Eu ^2$ or $\HH^2$, and let $n \geq 3$ be an odd integer.
If $C \subset \M$ is a disk and $P$ is a simple $n$-gon contained in $C$,
then there is an isosceles triangle, inscribed in $C$ and with side-lengths $\alpha \geq \beta$, such that $\perim P \leq (n-2)\alpha + 2\beta$.
\end{thm}

\begin{rem}\label{rem:spherical}
There is a real number $\varepsilon > 0$ such that for every odd integer $n \geq 3$ and $0 < \rho < \varepsilon$
the following holds. If $C \subset \S^2$ is a disk of radius $\rho$, and $P$ is a simple $n$-gon contained in $C$,
then there is an isosceles triangle, inscribed in $C$ and with side-lengths $\alpha \geq \beta$, such that
$\perim P \leq (n-2)\alpha + 2\beta$.
\end{rem}

In the hyperbolic and the spherical case, using simple calculus, it is possible to determine the maximum of the quantity $(n-2)\alpha + 2\beta$ as a function of $n$ and the radius of $C$ like in the Euclidean case.
Unfortunately, these functions seem too long and complex to be included in this paper.

\section{Proof of Theorem and Remark}

First, note that it is sufficient to find a triangle $T$, \emph{contained} in $C$, with side-lengths $\alpha, \beta$
and $\gamma$ such that $\perim P \leq (n-2)\alpha + \beta + \gamma$, since in this case we can move the vertices of $T$
to $\bd C$ in a way that no sidelength of $T$ decreases, and then we can move the vertex of $T$, opposite of the side with length $\alpha$, to the symmetric position while we do not decrease the quantity $(n-2)\alpha + \beta + \gamma$. 

We assume that if $\M = \S^2$, then the radius of $C$ is at most $\frac{\pi}{4}$.
Let $p,q \in C$ be distinct points, and if $\M = \S^2$, let $H$ be the open
hemisphere with the midpoint of $L(p,q) \cap C$ as its centre.
We introduce a partial order $\leq_{pq}$ on the points of $C$.
We set $p \leq_{pq} q$, and order the other points of $L(p,q)$ (or $L(p,q) \cap H$ in the spherical case) accordingly.
Let $x,y \in C$, and let $h_x$ and $h_y$ denote the orthogonal projections of $x$ and $y$ onto $L(p,q)$ (or $L(p,q) \cap H$ in the spherical case), respectively.
We say that $x \leq_{pq} y$ if, and only if, $h_x \leq_{pq} h_y$. If $h_x = h_y$ we say that $x =_{pq} y$.

\begin{defn}\label{defn:noncrossing}
Let $p,q,a,b,c \in \M$. We say that the segment $[p,q]$ and the polygonal curve $[a,b]\cup [b,c]$
\emph{do not cross} if the set $(p,q) \cap \big( (a,b] \cup [b,c) \big)$ is either empty or a segment.
\end{defn}

To prove Theorem and Remark, we need the following lemma.

\begin{lem}\label{lem:increasing}
Let $p,q \in C$ be distinct points and set $\delta = \dist_M (p,q)$.
Let $a,b,c \in C$ be points such that $a \leq_{pq} b \leq_{pq} c$,
and the polygonal curve $[a,b] \cup [b,c]$ and $[p,q]$ do not cross.
If $\dist_M (a,b) \leq \delta$ and $\dist_M (b,c) \leq \delta$, then there are points $a',b',c' \in C$
such that $\dist_M (a,b) + \dist_M (b,c) \leq \dist_M (a',b') + \dist_M (b',c')$, and $\delta \leq \dist_M (a',c')$.
\end{lem}

First, we show how our theorem and remark follow from Lemma.

Let $P$ be a simple $n$-gon contained in $C$. Let the vertices of $P$ be $a_1, a_2, \ldots, a_n=a_0$ in this order.
Without loss of generality, we may assume that $[a_0,a_1]$ is a longest side of $P$.
Set $p=a_0$ and $q=a_1$.

Consider the $(n+1)$-element sequence $\mu_1,\mu_2, \ldots \mu_{n+1}$ defined as follows.
For $i=1,2,\ldots, n+1$,
\[
\mu_i = \left\{
\begin{array}{l}
1, \quad \hbox{if } a_i \leq_{pq} a_{i+1} \hbox{ and } a_{i+1} \neq_{pq} a_i;\\
-1, \quad \hbox{if } a_{i+1} \leq_{pq} a_i \hbox{ and } a_i \neq_{pq} a_{i+1};\\
0, \quad \hbox{if } a_i =_{pq} a_{i+1}.
\end{array}
\right.
\]

Observe that $\mu_0 = \mu_n \neq 0$.
Since $n+1$ is even, there are two consecutive elements of the sequence that are both nonnegative or both nonpositive.
Thus, we may apply Lemma, and the assertion follows.

\section{Proof of Lemma}

First, we prove the lemma for $\M = \Eu^2$.
Consider a Descartes coordinate system and denote the coordinates of a point $z \in \Eu^2$
by $(\omega_z,\theta_z)$.
We may choose the coordinate system in a way that the centre of $C$ is the origin $o$, $\omega_p = \omega_q$, and
$\theta_p < \theta_q$.
Observe that for any $x,y \in C$, we have $x \leq_{pq} y$ if, and only if, 
$\theta_x \leq \theta_y$.
Let $R_p$ and $R_q$ be the two connected components of $L(p,q) \setminus (p,q)$ such that
$p \in R_p$ and $q \in R_q$.

Assume that $[a,b]$ intersects $R_q$.
Note that in this case $\dist_E(a,b) \leq \delta \leq \dist_E(b,p)$, and
$\delta \leq \dist_E(c,p)$.
Thus, we may choose $a'$, $b'$ and $c'$ as $p$, $b$ and $c$, respectively.
If $[b,c]$ intersects $R_p$, then a similar argument yields the assertion.

\begin{figure}[here]
\includegraphics[width=0.42\textwidth]{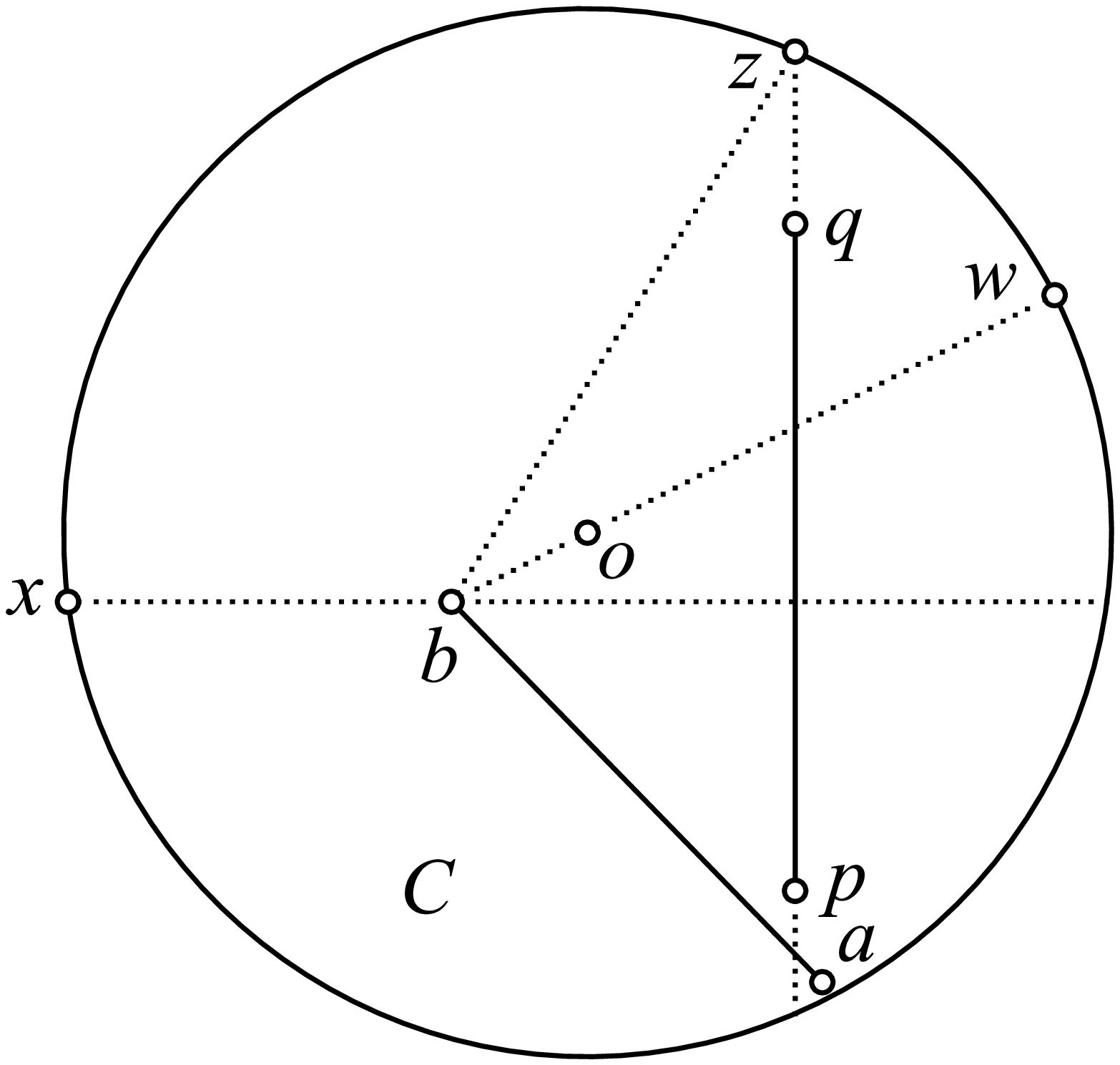}
\caption[]{}
\label{fig:Rp_intersect}
\end{figure}

Assume that $[a,b]$ intersects $R_p$.
If $[b,c]$ intersects $R_q$, then the assertion readily follows from $\dist_E(a,c) \geq \delta$.
Hence, we may assume that $b$ and $c$ are in the same closed half plane bounded by $L(p,q)$.
Let this half plane be denoted by $\bar{H}$, and let $\bar{C}$ denote the set of points with their $y$-coordinates
at least $\theta_b$ that are contained in $\bar{H} \cap C$.
Let $z$ be the intersection point of $R_q$ and $\bd C$, $w$ be the intersection point of
$\bd C$ and the ray emanating from $b$ and passing through $o$,
and $x$ be the point of $\bar{H} \cap \bd C$ with $\theta_x = \theta_b$ (cf. Figure~\ref{fig:Rp_intersect}).

If $w \in \bar{C}$, then $\dist_E(b,c) \leq \dist_E(b,w)$, and we may choose $a$, $b$ and $w$ as $a'$, $b'$ and $c'$, respectively.
Assume that $w \notin \bar{C}$, and observe that then $\dist_E(b,c) \leq \max \{\dist_E(b,x), \dist_E(b,z) \}$.
If $\dist_E(b,x) \leq \dist_E(b,z)$, then we may choose $a$, $b$ and $z$ as $a'$, $b'$ and $c'$, respectively.
If $\dist_E(b,z) \leq \dist_E(b,x)$, then we may choose $a$, $x$ and $z$ as $a'$, $b'$ and $c'$, respectively.
A similar argument proves the assertion in the case that $[b,c]$ intersects $R_q$.

We are left with the case that $a$, $b$ and $c$ are in the same closed half plane $H^+$ bounded by $L(p,q)$.
Without loss of generality, we may assume that $\area (H^+ \cap C) \geq \area (C \setminus H^+)$, and that $p, q \in \bd C$.
Now we drop the conditions that $\dist_E(a,b) \leq \delta$ and $\dist_E(b,c) \leq \delta$, and maximize
$\lambda$ under the conditions that $a,b,c \in C \cap H^+=C^+$
and $\theta_a \leq \theta_b \leq \theta_c$.

First of all, we may assume that $a, c \in C^+$.
Note that the tangent line of an ellipse is the external bisector of the angle bounded by the two rays, emanating from the point, that pass through the foci of the ellipse.
(This observation also holds in the hyperbolic plane and on the sphere.)
Thus, we may move also $b$ horizontally to $\bd C^+$ while not decreasing $\lambda$.
If $b \notin [p,q]$ and if $a$, $b$ and $c$ are in counterclockwise order in $\bd C^+$, then moving $a$ to $p$ and $c$ to $q$ does not decrease $\lambda$, and hence, the assertion follows.
If $b \notin [p,q]$ and if $a$, $b$ and $c$ are in clockwise order in $\bd C^+$, then we may move $b$ to $p$ or $q$ while not decreasing $\lambda$.

Assume that $b \in [p,q]$ and, without loss of generality, assume that $\theta_b \leq 0$.
Then $a$ may be moved to $H^+ \cap \bd C$ to the position such that $\theta_a = \theta_b$, and
$c$ may be moved to $H^+ \cap \bd C$ such that $o \in [b,c]$.

\begin{figure}[here]
\includegraphics[width=0.43\textwidth]{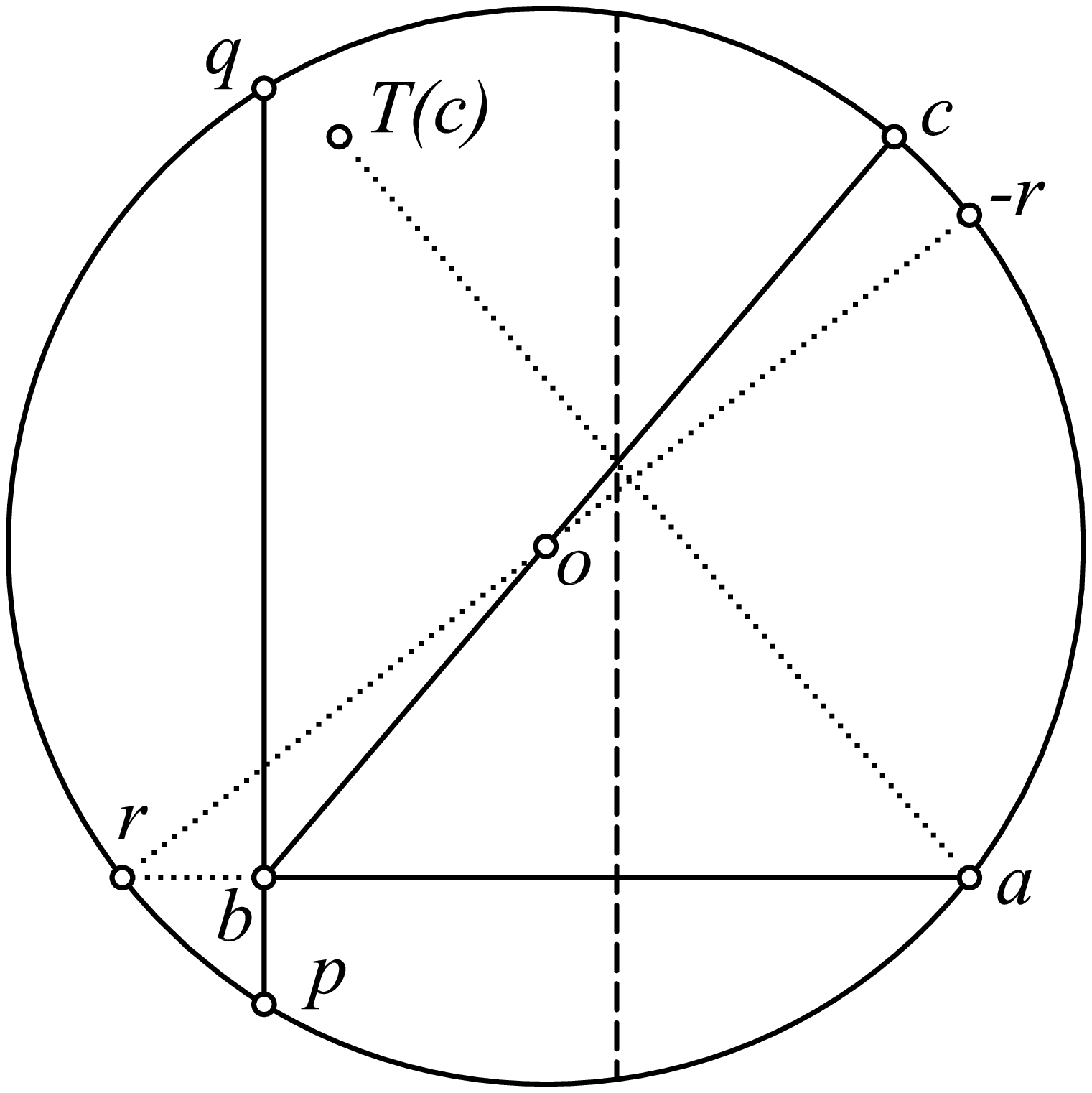}
\caption[]{}
\label{fig:samehalfplane}
\end{figure}

Let $r$ be the reflected image of $a$ about the $y$-axis, and observe that
$\omega_c \leq \omega_{-r}=\omega_a$ (cf. Figure~\ref{fig:samehalfplane}).
Let $T$ be the reflection about the line bisecting $[a,b]$.
Then $T(a)=b$, $T(b)=a$ and $T(c) \in C^+$.
Clearly, $\dist_E (a,b) \leq \dist_E(a,p)$ and $\dist(a,T(c)) \leq \dist_E(a,q)$, and
the assertion immediately follows.

Now we prove Lemma for $\M = \HH^2$.
The proof is a straightforward modification of the proof for $\M = \Eu^2$ in all the cases
but in the last one.
We prove the statement for this case only.
We define $H^+$ as in the Euclidean case.

Let $r$ be the intersection of $\bd C \setminus \inter H^+$ with the line $L(a,b)$,
and $s$ be the reflected image of $r$ about the centre of $C$.
Let $L_a$ be the line perpendicular to $L(a,b)$, that passes through $a$, and
let $t$ be the intersection point of $L_a$ and $\bd C$ different from $a$.
Note that $L_a$ does not separate $s$ and $[p,q]$, and that $L(p,q)$ is perpendicular to $L(a,b)$.
We observe that $q$, $c$, $s$ and $t$ are in this cyclic order in $\bd C$.

Let $T$ denote the reflection about the line bisecting the segment $[a,b]$,
and let $w$ be the intersection point of this line with $\bd C$ in the half plane, bounded by $L(p,q)$,
that contains $q$.
If $c$ is on the arc of $\bd C$, with endpoints $q$ and $w$, that does not contain $p$,
then from $\dist_H(a,b) \leq \dist_H(b,c) \leq \dist_H(b,a) + \dist_H(a,c)$ the assertion readily follows.
Assume the converse, that is, that $c$ is on the arc of $\bd C$, with endpoints $w$ and $t$, that does not
contain $p$.
Observe that $T(L(p,q))=L_a$, $T(L_a)=L(p,q)$, $T(a)=b$ and $T(b) = a$.
Note that the reflected image of the open arc in $\bd C$, with endpoints $t$ and $w$ and not containing
$p$ is an arc of radius $\rho$ that connects $w$ and a point of $[p,q]$.
Thus, $T(c) \subset H^+ \cap C$, and the assertion follows.

Finally, let $\M = \S^2$.
Let $R_p$ and $R_q$ be the two connected components of $(L(p,q) \cap H) \setminus [p,q]$
such that $p \in R_p$ and $q \in R_q$.

Assume that $[a,b]$ intersects $R_q$. Let this intersection point be denoted by $s$, and let $L_s$ be the line perpendicular to $R_q$ at $s$.
Let $t$ be the point of $H \cap L(p,q)$ such that $L_t=L(c,t)$ is perpendicular to $L(p,q)$ at $t$.
Observe that $t \in R_q$ and that it is not closer to $q$ than $s$.
As $\rho \leq \frac{\pi}{4}$, we obtain that the distance of $p$ from $L_s$ is equal to $\dist_S(p,s)$ and the distance of $p$ from $L_t$ is equal to $\dist_t(p,t)$.
Thus, $\dist_S(a,b) \leq \delta \leq \dist_S (p,b)$, and $\delta \leq \dist_S(p,c)$, and we may choose $p$, $b$ and $c$ as $a'$, $b'$ and $c'$, respectively.

The proof in the rest of the cases but the last one is a similarly modified version of the proof for $\M = \Eu^2$.
We prove the assertion for the last case.
We define $H^+$ as in the Euclidean case.

Let us imagine $\S^2$ as a unit sphere in $\Eu^3$.
Observe that for any $\nu > 0$, there is a value of $\varepsilon > 0$ such that if $\rho < \varepsilon$,
and $V$ is any stereographic projection of $\S^2$ onto a plane touching $\S^2$ at a point of $C$,
then, for any points $x,y \in C$, we have $| \dist_S(x,y) - \dist_E(V(x), V(y)) | < \rho \nu$.
We choose $\varepsilon$ to satisfy this condition, but we fix $\nu$ later.

Let $U$ denote the stereographic projection of $\S^2$ onto the plane of $\Eu^3$ that touches $\S^2$ at $b$.
Note that $U(C)$ is a Euclidean circle in $H$, and $L(U(a),U(b))$ is perpendicular to $L(U(p),U(q))$.
Furthermore, we have
\begin{equation}\label{eq:1}
\dist_S(a,b) + \dist_S(b,c) < \dist_E(U(a),U(b)) + \dist_E(U(b),U(c)) + 2\rho\nu.
\end{equation}

Let $m \in \bd C$ be the point such that the convex hull of $U(m)$, $U(p)$ and $U(q)$ is an isosceles triangle
that contains the centre of $U(C)$, and set $\beta = \dist_E(U(p),U(m))$.
Then
\begin{equation}\label{eq:2}
2\beta - 2\rho \nu < \dist_S(m,p) + \dist_S(m,q)
\end{equation}

Observe that there is a number $\tau > 0$ such that
\begin{equation}\label{eq:3}
\dist_E(U(a),U(b))+ \dist_E(U(b),U(c)) \leq 2 \beta - \rho \tau
\end{equation}
for all the possible positions of $a,b$ and $c$.
Thus, if $\nu \leq \frac{\tau}{4}$, then, from (\ref{eq:1}), (\ref{eq:2}) and (\ref{eq:3}),
we obtain that
\[
\dist_S(a,b) + \dist_S(b,c) < 2\beta + 2\rho \nu - \rho \tau \leq \dist_S(p,m) + \dist_S(m,q)
\] 
and the assertion follows.

\section{A remark and questions}

We remark that Lemma does not hold for spherical disks of any
radius $\rho < \frac{\pi}{2}$.
To show this, we have the following example.

\begin{example}
Let $\varepsilon > 0$ and let $p, q \in \S^2$ be two points with $\dist_S (p,q)=\pi - \varepsilon$. 
Let $k, -k$ be the centres of the two open hemispheres bounded by $L(p,q)$.
Let $[a,b]$ be a segment of length $\pi - \varepsilon$, perpendicular to $[p,q]$, such that $q$ is the midpoint of $[a,b]$. Note that $k,-k \in L(a,b) \setminus [a,b]$.
Thus, there is a disk $C$ of radius $\rho < \frac{\pi}{2}$ that contains $[p,q]$ and $[a,b]$.
Now we set $c=a$, and observe that $a \leq_{pq} b \leq_{pq} c$ and,
if $\varepsilon$ is sufficiently small, then $3\pi-3\varepsilon$ is clearly
greater than the perimeter of any triangle inscribed in $C$.
\end{example}

We ask the following questions.

\begin{ques}
Let $n \geq 5$ be odd, $0 < \rho < \frac{\pi}{2}$, and $C \subset \S^2$ be a disk
of radius $\rho$.
What is the supremum of the perimeters of the simple $n$-gons contained in $C$?
\end{ques}

It seems very difficult to find the maximal perimeters of simple $n$-gons contained in
an arbitrary given plane convex body.
Nevertheless, the following problem seems natural.

\begin{prob}\label{prob:oval}
Let $n \geq 5$ be odd, and let $C \subset \Eu^2$ be a plane convex body.
Prove or disprove that if $P$ is a simple $n$-gon contained in $C$, then there is
a triangle, inscribed in $C$ and with side-lengths $\alpha, \beta$ and $\gamma$, such that
$\perim P \leq (n-2)\alpha + \beta + \gamma$.
Is it true for plane convex bodies in the hyperbolic plane or on the sphere?
\end{prob}

We note that if $n$ is even, then the best possible upper bound for the perimeter of $P$
in Problem~\ref{prob:oval} is the trivial upper bound $n \diam C$.
This observation is valid also in any Minkowski plane.

\begin{ques}
Let $n \geq 5$ be odd, and let $\M$ be a Minkowski plane with the unit disk $C$.
What is the supremum of the perimeters of the simple $n$-gons contained in $C$?
In particular, can Theorem be generalized for Minkowski planes?
Can it be generalized for an arbitrary plane convex body of $\M$ instead of the unit disk of $\M$?
We note that, clearly, in these cases the optimal triangle inscribed in $C$ is not necessarily isosceles.
\end{ques}

{\bf Acknowledgements.}
The author is indebted to K. B\"or\"oczky jr., \'A. G. Horv\'ath and M. Nasz\'odi for their very helpful remarks.

\end{document}